%% file: paper.tex
\documentclass[11pt]{amsart}

\usepackage{amsfonts, enumerate} 
\usepackage{amssymb}
\usepackage[final]{graphicx, graphics,stmaryrd} 

\usepackage[margin=3cm]{geometry} 

\usepackage[breaklinks,bookmarksopen,bookmarksnumbered]{hyperref}

\usepackage[english]{babel}
\usepackage[utf8]{inputenc}
\usepackage{amsthm}
\usepackage{graphics}
\usepackage{amsmath}
\numberwithin{equation}{section} 
\usepackage{amstext}

\usepackage[safe,extra]{tipa}
\usepackage{multirow}
\usepackage{enumitem}

\usepackage{bm} 
\usepackage{stackrel} 
\usepackage{tikz-cd} 
\usepackage{mathtools} 

\input{newcommands_eng.tex}

\title{An extension to non-nilpotent groups of Rothschild-Stein lifting method}
 \author{Mattia Galeotti}
 \address{Mattia Galeotti, Dipartimento di Matematica,
         Università di Bologna\\
         Piazza di Porta San Donato, 5 - 40126 Bologna, Italy}
 \email{mattia.galeotti4@unibo.it}

\begin{document}

\begin{abstract}

 In their celebrated paper of 1976, Rothschild and Stein 
 prove a lifting procedure that locally reduces to a free nilpotent Lie algebra
 any family of smooth vector fields $X_1,\dots,X_q$, over a manifold $M$. 
 Then, a large class of differential operators can be lifted,
 and fundamental solutions on the lifted
 space can be re-projected to fundamental
 solutions of the given operators on $M$.
 In case that the Lie algebra $\g=\Lie(X_1,\dots,X_q)$ is finite dimensional but not nilpotent,
 this procedure could introduce a strong tilting of the space. 
 In this paper we represent a global construction of a Lie group $G$ associated to $\g$
 that avoid this tilting problem.
In particular $\Lie(G)\cong\g$ and a right $G$-action exists over $M$,
faithful and transitive, inducing a natural projection $E\colon G\to M$.
We represent
the group $G$ as a direct product $M\times G^z$ where the model fiber $G^z$ has a group structure.
We prove
that for any simply connected manifold $M$ -- and a vast class of non-simply connected manifolds --
 a fundamental solution for a differential operator $\mc L=\sum_{\al\in\N^q} r_\al\cdot X^\al$ of finite degree over $M$
 can be obtained, via a saturation method,
from a fundamental solution for the associated lifted operator
over the group $G$. This is a generalization
of Biagi and Bonfiglioli analogous result for homogeneous vector fields over $M=\R^n$.
\end{abstract}
\maketitle

{\tiny \gra{Keywords.} Partial Differential Equations; H\"ormander operators; Lifting and approximation theorem; Lie algebra of vector fields;
Lie group action;
Integration of vector fields.}

\section{Introduction}

Starting from the celebrated papers \cite{hor67, folla75, rothstein76},
a large literature has been devoted to the properties of 
the operator $\mc L=\sum_{i}^q X_i^2$
such that $X_1,\dots X_1$ are smooth vector fields over an orientable 
manifold $M$ and $\g=\Lie(X_1,\dots,X_q)$ verifies the H\"ormander's condition,
\cor{i.e.}~it spans at every point $x$ the tangent plane $T_xM$.

A key idea developed in these papers, is that the regularity properties of the solution are strictly related 
to the properties of the Lie algebra $\g$ generated by the vector fields. The H\"ormander's condition ensures 
the hipoellipticity of the operators, but more sophisticated properties of the solution depend on the stratification 
of $\g$.
In their crucial work \cite{rothstein76},
Rothschild and Stein introduced a ``lifting and approximation''
technique that allows to locally approximate
the operator $\mc L$ with a homogeneous left invariant operator, even
if we don't have a group structure on $M$.
Goodman develop a slightly different lifting in \cite{good78}.

Their method can handle any family of vector fields
satisfying the H\"ormander's condition and approximate it
with a nilpotent Lie algebra of vector fields.
This approximation is more than sufficient to  
find a fundamental solution or determine local regularity properties of the solution.
However, in this procedure global properties
such as symmetries, periodicity and invariance properties, are completely lost.
On the other side, preserving invariance properties with respect to suitable Lie groups (possibly non-nilpotent) 
is one of the main reasons for studying differential equation in Lie group structures. 
Hence a challenging problem is to associate to the given family of vector fields, new models
which provide global approximation of the structure, global properties of the fundamental solution, 
and preserve the invariance properties.\newline

Let us consider a simple example.
Let $X_1,X_2$ be smooth vector fields in $M=\R^2$
such that $X_1,X_2$ and $X_3=[X_1,X_2]$ span the whole space at every point.
Also assume that $[X_1,X_3]=-X_2$ and $[X_2,X_3]$. This set of relations
allows to define a Lie algebra structure $\mfk g$ over $\langle X_1,X_2,X_3\rangle$.
In this relations it is hidden a periodicity condition encoding
meaningful properties of the solution of $\mc L=X_1^2+X_2^2$.
As a consequence, these vector fields admit a natural representation in terms of the 
sinus and cosinus functions, or equivalently in terms of the simply connected Lie group
of rigid motions, which is the universal cover $G\cong \R^3$ of $SE(2)\cong \R^2\times S^1$.
See also Example~\ref{exsin} and observe that $S$ is the unique connected and simply connected
Lie group such that $\Lie(G)=\g$.

However, in the Rothschild-Stein and Goodman procedures, this properties are lost, 
and the coefficients of the local approximations of $X_1,X_2,X_3$ are polynomials,
 the Taylor developments of the sinus or
cosinus function around a chosen point.
In this work
we introduce a new framework for a global action on $M$ of the Lie group $G$ associated to $\g$,
without any nilpotency condition. In this way, we ``lift''
the study of the operator $\mc L$ to a left invariant operator $\widetilde {\mc L}$ on $G$,
obtaining the solution of the first as a ``saturation'' of the solution of the latter.\newline

In full generality, our main objects of study are differential operators $\mc L$ such that 
\begin{equation}\label{diffop}
\mc L(x)=\sum_{|\al|\leq k}r_\al(x)\cdot X^\al(x),\ \ \forall x\in M,
\end{equation}
where any $\al\in \N^q$ is a multi-index of length bounded by $k>0$
and the $r_\al\colon M\to \R$ are smooth coefficients.
If a homogeneous group structure exists over $M=\R^m$
such that $\mc L$ is left invariant and homogeneous of degree $2$, then
Folland built in \cite{folla75} a homogeneous fundamental solution for~$\mc L$.
In \cite{rothstein76}, Rothschild and Stein 
build
a higher dimensional nilpotent group $G_U$ for some
bounded open subsets $U\subset M$ of a covering family
for the variety $M$,
and for any such  neighborhood, a natural projection $E\colon G_U\to U$ exists.
The idea is that we can lift the vector fields $X_i$ to $G_U$
and approximate them with elements of $\Lie(G_U)$ (see also \cite{mitch85}).
These tools allow estimates on the Sobolev norm of $u$ and $\mc Lu$ for any distribution
$u$ over $U$, and the study of the hypoellipticity of $\mc L$.
For a focus on the Sobolev estimates, see \cite{brabra23}.
In a similar fashion,
in the Goodman's procedure
 the group $G_U$ is generated by vector fields $\widetilde X_i$ that
are $E$-related to the $X_i$, meaning that
$\widetilde X_i(f\circ E)=X_if\ \ \forall f\in C^\infty(M)$.\newline

Rothschild-Stein and Goodman are both local approaches, meaning that any fundamental solution
for $\mc L$ is built over the separate neighborhoods $U$, and then glued together.
The works
\cite{folla75, folla77} by Folland are instead focused on a \cor{global} approach on $M=\R^m$
in the case that the vector fields $X_i$ are homogeneous with respect
to a set of dilations $\delta_\lambda$.
In this case  a global lifting
$E\colon\mathbb G\to \R^m$ exists, where $\mathbb G$ is a Carnot group
generated by the vector fields $\widetilde X_i$ as in the Goodman's lifting,
but over the whole variety.
As a consequence, the Carnot sub-Laplacian $\mc L_{\mathbb G}$ and the operator $\mc L=\sum X_i^2$
are $E$-related (for a wider study of this operator see also the monograph \cite{blu07}).
The advantage of the Folland's approach is that it allows a global representation
and therefore it improves the knowledge of the associated Sobolev spaces.


As Folland and Stein treated in \cite{follastein74, folla75}, over a Carnot group the sub-Laplacian $\mc L_{\mathbb G}$
has always a fundamental solution $\Gamma_{\mathbb G}$.
The question then arises, about the existence of a fundamental
solution $\Gamma$ for $\mc L$ which is in some sense $E$-related to $\Gamma_{\mathbb G}$.
Biagi and Bonfiglioli answer positively to this question in \cite{biabon17},
by a saturation process, meaning an integration along the $E$-fibers.
In their work they consider any differential operator in the form \eqref{diffop} over $\R^m$.
 The Carnot group $\mathbb G$ is then represented through a change of coordinates
as a product $\R^m\times \R^p$, and the projection $E\colon \mathbb G\to \R^m$ becomes
the first coordinate projection. If $\widetilde {\mc L}$ is the lifting of $\mc L$ and $\widetilde \Gamma$
the fundamental solution for $\widetilde{\mc L}$, then the ``saturated'' function
\begin{equation}\label{satu}
\Gamma(x;y)=\int_{\R^p}\widetilde\Gamma(x,0;y,t)\d t\ \ \forall x,y\in \R^m,\ x\neq y
\end{equation}
is in fact the fundamental solution for $\mc L$ over $\R^m$.
 Their construction also allows for estimates on the original operator~$\mc L$,
developed for example in \cite{bbb21, bbb22}.\newline

In the present work we continue the program initiated by Biagi and Bonfiglioli
by considering a large class of differential operators over a manifold $M$,
and dropping the homogeneity conditions on the vector fields $X_i$.
In particular, the differential operators considered are in the form \eqref{diffop}
with the Lie algebra $\Lie(X_1,\dots,X_q)$ that is complete and finite dimensional
but not nilpotent. 
This allows for global
representations of the fundamental solutions even without a dilation structure on $M$. 

The central tool of our analysis
is the Fundamental Theorem of Lie Algebra Actions (see Theorem \ref{fund}),
whose usefulness in this setting has been emphasized by the same Biagi and Bonfiglioli in \cite{biabon23}.
If $G$ is any finite dimensional Lie group with a right action $\mu\colon M\times G\to M$,
then its differential induces a Lie algebra morphism $T_\mu\colon \Lie(G)\to \mfk X(M)$
called infinitesimal generator of $\mu$.
The Fundamental Theorem states the converse,
if 
$T\colon \Lie(G)\to \mfk X(M)$ is a complete Lie algebra morphism, then there exists
a unique right $G$-action $\mu\colon M\times G\to M$ whose infinitesimal generator is $T$.
If $G$ is the unique connected and simply connected Lie group $G$ such that 
$\Lie(G)=\Lie(X_1,\dots,X_q)$,
then the H\"ormander's condition corresponds to the transitivity of the induced
right action $\mu\colon M\times G\to M$. By fixing a starting point $z\in M$, we obtain
the Folland-style morphism $E:=\mu(z,-)\colon G\to M$. 
We are now able to show a product representation
\[
G\cong M\times G^z,
\]
where $G^z=E^{-1}(z)$ is the central $G$ fiber and 
has a group structure.
This construction allows for a general saturation procedure and therefore
 the main
result of our work (see Theorem \ref{main} for the details).

\begin{teon}
 If $\mc L$ is a differential operator
of the kind \eqref{diffop}, $\widetilde {\mc L}$ its lifting through $E$ and 
$\widetilde \Gamma$ is a fundamental solution for $\widetilde{\mc L}$, then
\begin{equation}\label{satu2}
\Gamma(x;y)=\overline \rho(x)\cdot \int_{G^z}\widetilde\Gamma(x,s;y,t)\d t, \ \forall x,y\in M\ x\neq y,
\end{equation}
is a fundamental solution for $\mc L$, where $\overline \rho\colon M\to \R$ is a smooth
everywhere non-null function and $s$ any element of the group $G^z$.
\end{teon}
It has to be noted that, by exploiting the action of the central fiber $G^z$
on any other $E$-fiber, we can prove that 
\[
\widetilde{\mc L}^*=\mc L^*+\sum_{|\gamma|\geq 1}r_{\beta,\gamma}^*(x)\cdot X^\beta Y^\gamma,
\]
where the vectors $Y$ are defined along the $G^z$ coordinate, and the coefficients depend only on the $M$ coordinate.
Therefore we achieve a separation of the coordinates such that the representation
of  the differential operator is invariant along the fibers
of the $E$ map.
This ``vertical invariance'' of the lifting is a key result, allowing the saturation procedure
even without a dilation structure on $M$ or $G$.\newline

In Section \ref{basic} we introduce the notation and the Fundamental theorem of Lie algebra actions.
In \ref{coo} the group $G$ generated by the given vector fields on $M$,
is represented as a direct product between $M$ and a group fiber.
In Section \ref{sfs} the existence theorem for the ``saturated'' fundamental solution
 is stated and proved for simply connected manifolds, while in Section
\ref{nonsimply} we consider a class of non-simply connected manifolds where
the same result is still valid.\newline

\gra{Aknowledgements.}
The author wants to thank Giovanna Citti for her help in this and other research endeavors, and Andrea Bonfiglioli
for the crucial suggestions in perfecting this work.

\section{Basics}\label{basic}

\subsection{Basic notations}
We consider a smooth simply connected manifold $M$
and a finite family of vector fields over it. We denote
by $\mfk X(M)$ the space of smooth vector fields over $M$.
Given $X_1,X_2,\dots, X_q\in \mfk X(M)$,
if $\al=(\al_1,\dots,\al_q)$ is a multi-index in $\N^q$,
we use $X^\al$ as a short notation for
\[
X^\al=X_1^{\al_1}X_2^{\al_2}\cdots X_q^{\al_q}.
\]
Let $\g:=\Lie(X_1,\dots,X_q)$ be the Lie algebra generated by the $X_i$s,
that is the algebra generated by the $X_i$s
and their commutators of any order
\[
[X_{i_1},[X_{i_2}[\dots,X_{i_k}]]],
\]
then we suppose that $\g$ respects three main conditions stated here.

\begin{rmk}\label{rmk1}
These conditions are in fact a completeness condition and a re-writing of the well known H\"ormander's condition:
\begin{enumerate}
\item $\g$ is finite dimensional, if necessary we will complete the sequence
above to a basis $\{X_1,\dots,X_n\}$;
\item for any $X\in \g$, $X$ is a complete vector field;
\item at every $x\in M$, $\g_x:=\{X_x|\ X\in \g\}$ is the whole
tangent space $T_x M$.
\end{enumerate}
\end{rmk}
If the three conditions are respected and we denote by $m$ the $M$ dimension,
then by definition~$m\leq n$.\newline

We use the following notation for a differential operator
$\mc L$ of order $k$ over $M$,
\[
\mc L=\sum_{|\al|\leq k} r_\al\cdot X^\al,
\]
where $|\al|$ is the length of the multi-index $\al$,
and the $r_\al$ are real smooth coefficients.\newline

By Lie's Third Theorem (see \cite{ser09}) we know that there exists unique
a connected and simply connected Lie group $G$ such that 
$\Lie(G)\cong \g=\Lie(X_1,\dots,X_q)$.
We denote by $L\colon \Lie(G)\to \g$ the isomorphism,
and sometimes we will use the notation $\widetilde X$ for $L^{-1}X$ 
with $X$ any vector field in $\g$. Therefore $L\widetilde X=X$.\newline

\subsection{Lie algebra maps and right group actions}
The key ingredient of our work is the Fundamental Theorem of Lie Algebra Actions. We recall it here. 
Consider a simply connected Lie group $G$, a smooth manifold $M$
and a right $G$-action 
\[
\mu \colon M\times G\to M.
\]
When it is clear from the context, we will denote by $z\cdot \xi$ the action $\mu(z,\xi)$.
We denote by $\mu^{(z)}$ the map $\mu(z,-)\colon G\to M$.

This action induces a map from the Lie algebra associated to $G$
to the vector space $T_zM$ for any point $z$ of $M$.
Indeed, to $\widetilde X\in \Lie(G)$ we can associate the vector
\[
\left.\frac{\d}{\d t}\right|_{t=0}\mu(z,\Exp_G(t\widetilde X))=\d_e \mu^{(z)}(\widetilde X).
\]
Moreover we observe that, if we denote by $L_\xi\colon G\to G$ the left $\xi$ multiplication
inside $G$, then
\[
\mu^{(z\cdot \xi)}=\mu^{(z)}\circ L_\xi,
\]
this implies that 
\[
\d_e\mu^{(z\cdot \xi)}=\d_\xi\mu^{(z)}\circ\d_eL_\xi\colon \Lie(G)\to T_{z\cdot \xi}M.
\]
As any vector $\widetilde X$ in $\Lie(G)$ is left-invariant,
the formula above shows that
\begin{equation}\label{dthet}
\d_e\mu^{(z\cdot \xi)}\widetilde X=\d_\xi\mu^{(z)}\widetilde X.
\end{equation}
This allows to define a map 
\[
T_\mu\colon \Lie(G)\to \mfk X(M)
\]
which is called the \cor{infinitesimal generator} of the $\mu$ action.
The important result about the infinitesimal generator $T_\mu$ is the fact that it is a Lie algebra morphism (\cite[Theorem 20.15]{lee12}).
Given any Lie algebra $\h$ and a Lie algebra morphism $T\colon \h\to\mfk X(M)$,
this morphism is called complete if for any $\widetilde X\in \h$, the vector field $T(\widetilde X)$
is complete.

\begin{teo}[Fundamental Theorem of Lie Algebra Actions, {\cite[Theorem 20.16]{lee12}}]\label{fund}
Let $M$ be a smooth manifold, $G$ a simply connected Lie group and $T\colon \Lie(G)\to \mfk X(M)$
a complete Lie algebra morphism. Then, there exists a unique right $G$-action $\mu\colon M\times G\to M$
whose infinitesimal generator is $T$.
\end{teo}
Thanks to the above theorem, if we consider the aforementioned isomorphism
\[
L\colon \Lie(G)\to\g\subset\mfk X(M),
\]
this is the infinitesimal generator of a right $G$-action
\[
\mu_L\colon M\times G\to M.
\]
As already said, we denote by $E\colon G\to M$ the map
$\mu_L^{(z)}$. As a consequence of \eqref{dthet} and the definition
of the $L$ isomorphism, for any $X\in \g$,
$\widetilde X=L^{-1}X$ and $X$ are $E$-related, meaning that
\begin{equation}\label{erel}
\d_\xi E(\widetilde X)=X_{E(\xi)}.
\end{equation}
This is also proved by Biagi and Bonfiglioli, \cite[Formula (3.10)]{biabon23}.

Observe that, by construction, the $E$ map
respects the equality
\begin{equation}\label{eqpsi}
E(\xi)=\Psi_1^{L\log_G(\xi)}(z)
\end{equation}
for $\xi$ in an opportune neighborhood of the $G$ identity element $e$,
where $\Psi_t^X$ is the usual flow operator along the vector field $X$.\newline

 We denote by $G^x:=E^{-1}(x)$
the fiber over any point $x\in M$ and observe that in the case of $z=E(e)$, the fiber $G^z$ is a subgroup of $G$.
Moreover, it is possibile to define a $G^z$-action on any fiber $G^x$ by left multiplication. Indeed, if $\xi\in G^z$ and $\eta\in G^x$,
$\xi\eta\in G^x$ for any $x\in M$. This action is faithful and transitive.

\begin{rmk}\label{sc}
If $M$ is simply connected, then the $E$-fibers are connected too.
Indeed, if a path $\gamma$ relying two connected components
of $G^x$ exists for some $x\in M$, then its image $E(\gamma)$ would
be an un-contractible loop.
\end{rmk}
\begin{ex}\label{exgru}
Consider the case of the Grushin operator $\mc L=\de_{x_1}^2+\left(x_1\de_{x_2}\right)^2$
defined over~$\R^2$.
We consider the smooth and complete vector fields $X_1=\de_{x_1}$ and $X_2=x_1\de_{x_2}$. Even
if they don't generate linearly the tangent space at $(0,0)$, their Lie algebra respects the 
H\"ormander's condition everywhere.
Indeed, $\g=\Lie(X_1,X_2)=\langle X_1,X_2,X_3\rangle$
with $X_3=[X_1,X_2]=\de_{x_2}$, which span $\R^2$ at every $(x_1,x_2)\in \R^2$.

The connected and simply connected group whose Lie algebra is $\g$
is the first Heisenberg group $\H^1$. If we start at $z=(0,0)$ 
setting $\xi_1,\xi_2,\xi_3$ as exponential coordinates of the first type on $\H^1$
with respect to the basis $X_1,X_2,X_3$, then
\[
E(\xi_1,\xi_2,\xi_3)=\exp_{(0,0)}(\xi_1X_1+\xi_2X_2+\xi_3X_3)=\left(\xi_1,\xi_3+\frac{\xi_1\xi_2}{2}\right),
\]
where $\exp$ is the usual exponential map on $\R^2$.\newline
\end{ex}

\section{Coordinate system on the group}\label{coo}
\subsection{The group as a fiber bundle}\label{coorP}
Given the isomorphism $L\colon \Lie(G)\to \g$ as above,
a simply connected smooth manifold $M$
and a point $z\in M$,
we consider the map $E=\mu_L^{(z)}$
and the $z$ stabilizer $G^z$.

\begin{rmk}
By the hypothesis in Remark \ref{rmk1},
for any $z\in M$ there exists a vector subspace $\g'\subset \g$
of dimension $m$ and a neighborhood $U$ of $z$ such that
$\g'_x=T_xM$ for any $x\in U$. Therefore it is an exponential-kind
diffeomorphic (in a neighborhood of the origin) map
\begin{align*}
\g'&\to M\\
X&\mapsto \Psi_1^X(z).
\end{align*}
We denote by $\overline\log\colon U\to \g'$ its inverse. This allows to define a trivialization of the
map $E\colon G\to M$. Indeed, if $G_U:=E^{-1}(U)$, consider the map $G_U\to U\times G^z$,
\[
\xi\mapsto \left(E(\xi),\ \xi\cdot \Exp_G(-L^{-1}\overline\log (E(\xi)))\right).
\]
This is well defined as a consequence of \eqref{eqpsi}.
Equivalently, for any $\xi\in G$ we can write it as 
$\xi=s\cdot \Exp_G(\widetilde X)$ where $s\in G^z$ and $\widetilde X$ is in $\g'$.

We built a local trivialization of the bundle $G\to M$, proving it is a fiber bundle.
\end{rmk}

\begin{rmk}
The \cor{vertical} bundle associated to $E$ is the subbundle of $TG\to M$ defined
as $V(E):=\Ker(\d E)\subset TM$.
An Ehresmann connection on a fiber bundle is the data of a \cor{horizontal bundle}, meaning
a sub-bundle $H(E)\subset TM$ such that $V(E)\oplus H(E)=TM$. Observe that this implies $H_\xi(E)\cong T_{E(\xi)}M$
for any $\xi\in G$.

The differential $\d E$ is surjective at any $\xi\in G$, therefore if we have an Ehresmann connection $H(E)$ and
if $x=E(\xi)$, then there exists a neighborhood $U\subset M$ of $x$ and a smooth subvariety $\widetilde U \subset G$
containing $\xi$
such that $E|_{\widetilde U}$ is a diffeomorphism between $\widetilde U$ and $U$,
and $T\widetilde U=H(E)|_{\widetilde U}$.
Over $U$ we can therefore define $\ell\colon U\to G$ as the $E$ inverse.
By simple connectedness of $M$, we can impose $\ell(z)=e$ the identity element of $G$,
and extend $\ell$ to a section $\ell\colon M\to G$
on the whole manifold.
\end{rmk}
\begin{ex}
We consider again the case of the Grushin operator and the induced group $\H^1$ acting on $\R^2$,
as in Example \ref{exgru}.
As showed, the map $E\colon \H^1\to \R^2$ gives two coordinates
\[
\left\{\begin{array}{ll}x_1&= \xi_1\\
x_2&=\xi_3+\frac{\xi_1\xi_2}{2}\\
\end{array}\right.
\]
With $X_1,X_2,X_3$ as before, we denote by $\widetilde X_1,\widetilde X_2,\widetilde X_3$ their
abstract counterparts, meaning that we are looking at the same vector fields as elements of $\Lie(G)$
instead of $\g\subset \mfk X(\R^2)$.

In this case $V(E)=\langle \widetilde Y=\widetilde X_2-x_1\widetilde X_3\rangle$ by construction,
and the Biagi-Bonfiglioli construction is based on the Ehresmann connection $H(E)=\langle \widetilde X_1,\widetilde X_3\rangle$.
If we denote by $x_3$ the coordinate associated to $\widetilde Y$, this is the coordinate system
treated also in \cite{biabon17}, and in this representation $\H^1\cong \R^3$ with $E$ corresponding to the coordinate projection.
Observe that the liftings of $X_1,X_2$ are
\begin{align*}
L^{-1}X_1&=\de_{x_1}\\
L^{-1}X_2&=x_1\de_{x_2}+\de_{x_3},
\end{align*}
therefore the two lifted vector fields are everywhere linearly independent.\newline
\end{ex}

The section $\ell$ above allows to build a trivialization of $G$. 
Observe that $G^x=G^z\ell(x)$ by construction, therefore we have
\begin{align*}
\psi\colon G& \xrightarrow{E\times R_{\ell\circ E}^{-1}} M\times G^z\\
\xi&\mapsto (E(\xi), \xi\cdot \ell(E(\xi))^{-1})
\end{align*}
as a smooth isomorphism. 
Therefore the inverse diffeomorphism defines a different coordinate system on $G$,
\begin{equation}\label{csG}
(x,s)\mapsto s\cdot \ell(x)\in G.
\end{equation}
In this setting, the fiber $G^x$
is identified with $\{x\}\times G^z\subset M\times G^z$. Moreover,
if $E(\xi)=x$, then
 the following equalities hold for the tangent space to the $G^x$-fiber,
\begin{equation}\label{Gy}
T_\xi G^x=\d R_{\ell(x)}\left(T_{\xi\ell(x)^{-1}}G^z\right)=\Ker(\d_\xi E).
\end{equation}
We thus represented $G$ as a smooth fiber bundle over $M$ such that every fiber
is identified with the $G^z$ group.\newline

Let
$\omega_G$ be the left invariant volume form over $G$ (unique up to a multiplicative constant).
 If $\widetilde X_1,\dots,\widetilde X_n$ is a basis of $\Lie(G)$
and $\xi_1,\dots,\xi_n$ the exponential coordinates of the first kind induced on $G$ by the exponential
map with respect to this basis,
then we can set
\[
\omega_G=\d \xi_1\wedge\cdots \wedge \d  \xi_n.
\]
Moreover,
we consider the scalar product that makes $\widetilde X_1,\dots,\widetilde X_n$ in an orthonormal basis.
At every point $\xi\in G$, this scalar product induces a split
\begin{equation}\label{split2}
\Lie(G)=\Ker(\d_\xi E)^\bot\oplus\Ker(\d_\xi E).
\end{equation}
Then, considering the definition of $\psi$ and the identification \eqref{Gy}, there is a natural identification 
\begin{equation}\label{natid}
T_\xi G=\Ker(\d_\xi E)^\bot\oplus T_\xi G^x\xrightarrow{\d \psi=\d E\oplus \d R_{\ell(x)}^{-1}} T_xM\oplus T_{\xi\ell(x)^{-1}}G^z=T_{(x,\xi\ell(x)^{-1})}(M\times G^z).
\end{equation}
If for every $\xi\in G$, $\omega_{\Ker}(\xi)$ and $\omega_\bot(\xi)$ are the naturally induced volume forms
on $\Ker(\d_\xi E)$ and $\Ker(\d_\xi E)^\bot$ respectively, then the following wedge formula is also true,
\[
\omega_G=\omega_\bot\wedge\omega_{\Ker}.
\]
\begin{rmk}
Observe that once we fix the orientation of $\widetilde X_1,\dots,\widetilde X_n$,
$\omega_G$ is univocally determined while $\omega_\bot$ and $\omega_{\Ker}$
are determined together up to orientation
\end{rmk}

In the following sections we are going to 
precise how the coordinate change to $M\times G^z$ acts
on the form $\vol_M\wedge\omega_{\Ker}$.
In particular, we are going to prove the following important coordinate change formula.
\begin{teo}\label{teoimp}
There exists a smooth everywhere non-null function $\overline\rho\colon M\to \R$,
such that at any point $\xi\in G$ with $x=E(\xi)$,
\[
\omega_G(\xi)=\overline\rho(x)\cdot \d_\xi \psi^*\left(\vol_M\wedge\omega_{\Ker}\right).
\]\newline
\end{teo}

\subsection{Coordinate change on the base manifold}
We revisit the coordinate change rule on volume forms in order
to adapt it to our situation.

If $V$ is a (real) vector space of finite dimension $m$,
for any basis $\bm v=v_1,\dots,v_m$ of $V$, we denote by $\omega_{\bm v}$ the $m$-form
such that $\omega_{\bm v}(v_1,\dots,v_m)=1$.
Consider a scalar product $\langle,\rangle$ over~$V$, if $u\in V$
 we denote by $u^*$ the dual linear map $\langle u,-\rangle\in V^*$.
For example if $\bm e= e_1,\dots,e_m$
is an orthonormal basis of $V$, then
\[
\omega_{\bm e} =e_1^*\wedge\cdots\wedge e_m^*.
\]
This form is called volume form of $V$ and sometimes also denoted by $\omega_V=\omega_{\bm e}$.
By fixing the sign of the volume form,
we chose the orientation of any orthonormal basis.
If $V$ is the tangent space to a Riemannian oriented manifold $M$ and
the scalar product on it is the associated Riemannian metric,
then we denote the volume form by $\vol_M$.\newline

If $B\colon V\to V$ is a linear isomorphism, we use the notation
$B^*\omega$ or $\omega(B-)$ for the form 
\[
B^*\omega(w_1,\dots,w_m)=\omega(B(w_1,\dots,w_m))=\omega(Bw_1,\dots,Bw_m).
\]
As $\omega(B-)=\det(B)\cdot\omega$, if $B\bm v$ is the basis $Bv_1,\dots,Bv_m$,
we have $B^*\omega_{B\bm v}=\omega_{B\bm v}(B-)=\omega_{\bm v}$ and therefore
\begin{equation}\label{eq_bv}
\omega_{B\bm v}=\frac{1}{\det B}\cdot \omega_{\bm v},
\end{equation}
which is the classic coordinate change formula.
\begin{rmk}
In this setting we recall a well known volume formula.
Indeed, given a basis $\bm v=v_1,\dots, v_m$ of $V$, if $B$ is the map
sending $e_i$ to $v_i$ for any $i=1,\dots, m$, then
$\omega_{\bm v}=\omega_V(B^{-1}-)$, which gives
 $\omega_V=\det(B)\cdot \omega_{\bm v}$.
We observe that $\det(B)=\sqrt{\det(BB^T)}$
and $BB^T$ is precisely the matrix of the scalar products $g_{ii'}=(\langle v_i,v_{i'}\rangle)_{ii'}$,
thus
\[
\omega_V=\sqrt{\det(g_{ii'})}\cdot \omega_{\bm v},
\]
where we supposed that $v_1,\dots, v_m$ is positively oriented
and therefore $\det(B)=\det(B^T)>0$.\newline
\end{rmk}

We show a slightly more general result
of the one in the last remark.
If $V'$ is another 
(real, finite dimensional) vector space of dimension $m'$ equipped with a scalar product
and $B\colon V'\to V$ a linear surjective map, let $\omega_\bot$ be the
volume form on $\Ker(B)^\bot$.
\begin{lemma}\label{lem_bbt}
If $\omega_V$ is the volume form on $V$, then
\[
B^*\omega_{V}=\pm\sqrt{\det(BB^T)}\cdot \omega_{\bot}.
\]
\begin{rmk}
Here the $\pm$ sign means that the orientation 
determined by $B$ is not canonical.
\end{rmk}
\end{lemma}
\proof
Consider an orthonormal positively oriented basis $\bm v$ of $\Ker(B)^\bot$,
then by construction $B\bm v$ is a basis of $L$,
and $\omega_{B\bm v}(B-)=\omega_{\bm v}=\omega_{\bot}$.
Moreover, consider the map $C\colon V\to \Ker(B)^\bot$ sending the orthonormal
basis $\bm e$ of $L$ in $\bm v$,
then $\omega_{\bm e}=\omega_V$ and
\[
\omega_{B\bm v}=\omega_{BC\bm e}=\frac{1}{\det (BC)}\omega_{V}.
\]
By construction the matrix $C$ is orthogonal and
$\det(BB^T)=\det(BCC^TB^T)=\det(BC)^2$,
which proves the lemma.\fine

We are going to apply Lemma \ref{lem_bbt}
in our case, with the linear map given by
the restriction of the differential $\d_\xi E\colon \Ker(\d_\xi E)^\bot\to T_xM$,
where $x=E(\xi)$ as usual.

\begin{rmk}\label{rmkimp}
Observe that 
by the Formula \eqref{erel}
the image of $\d_\xi E$ only depends on the chosen vectors
on $\Lie(G)$ and the point $E(\xi)$, \cor{i.e.}~it
is a function of $E(\xi)$ and not $\xi$.
In particular,
\begin{equation}
\d_\xi E=\ev_{E(\xi )}\circ L.
\end{equation}
In the following we use the notation $\mc X(x):=\d_\xi E\colon \Lie(G)\to T_xM$
where $x=E(\xi)$.
\end{rmk}
\begin{rmk}\label{rmkimp2}
This point is crucial for the following of our work, in particular because
it states that the subspaces $\Ker(\d_\xi E)\subset \Lie(G)$ is invariant along
the $G^x$ fibers. Moreover, for any $x\in M$ we have the isomorphism 
\begin{equation}
\Ker(\mc X(x))=\Ker(\d_\xi E)\cong \Lie(G^z)\subset \Lie(G),
\end{equation}
but for different $x\in M$ it corresponds to different inclusions of $\Lie(G^z)$ inside $\Lie(G)$.\newline
\end{rmk}

If  $X_1,\dots, X_n$ is a basis of $\g$,
and $\widetilde X_1,\dots, \widetilde X_n$ the ``lifted''
basis on $\Lie(G)$,
then the map $\mc X$ is represented
by the $m\times n$ matrix 
\[
\mc X(x)=(X_1(x)|\dots|X_n(x)).
\]
Moreover, $\mc X^T(x)$ is the adjoint operator, and therefore
\[
\Ker(\d_\xi E)^\bot=\Ker(\mc X(x))^\bot=\im(\mc X(x)^T).
\]
We can conclude by Lemma \ref{lem_bbt}, obtaining
the coordinate change formula for the $M$ component in the isomorphism
$\psi\colon G\xrightarrow{\sim} M\times G^z$. 
\begin{equation}\label{eqE}
\left(\left.\d_\xi E\right|_{\Ker(\d E)^\bot}\right)^*\vol_M=\pm\sqrt{\det(\mc X(x)\mc X^T(x))}\cdot \omega_\bot.
\end{equation}
As already said, the $\pm$ sign means that the orientation induced by $\mc X(x)$
is not canonical. Anyway, up to re-ordering the basis $X_1,\dots, X_n$ of $\g$,
we can assume that the sign above is $+$ in an opportune neighborhood of $M$,
and therefore by orientability of $M$ the scalar factor is positive for any point of~$M$.
\begin{rmk}\label{rmkH}
The determinant above is everywhere non-null because of the H\"ormander-like hypothesis
on the vector fields $X_1,\dots,X_n$.\newline
\end{rmk}

\subsection{Coordinate change on the fibers}
As seen in \eqref{Gy},  the right
multiplication $R_{\ell(x)}$ induces an isomorphism
between the tangent spaces of the fibers $G^z$ and $G^x$.
Therefore, if we consider the action induced on the volume form of
this tangent spaces, that we called $\omega_{\Ker}$, we obtain
\[
\d R^*_{\ell(x)}\left(\omega_{\Ker}(\xi)\right)=c(\xi)\cdot \left(\omega_{\Ker}(\xi\ell(x)^{-1})\right),
\]
where $x=E(\xi)$ and $c$ is an everywhere non-null smooth function.
\begin{lemma}\label{lemc}
The scaling function $c$ is constant on any fiber $G^x$ 
of $G\xrightarrow{E} M$.
\end{lemma}
\proof
By construction the scalar product on $G$, and therefore the volume form, 
is invariant by left multiplication. As pointed out in Remark \ref{rmkimp2}
the subspace $\Ker(\d E)\subset \Lie(G)$ is invariant along the fibers $G^x$,
meaning it is invariant under left multiplication by elements in the subgroup $G^z\subset G$.

All this implies that the volume form $\omega_{\Ker}$
is invariant by left multiplication by $G^z$ too.
That is, for any $s\in G^z$,
\[
\omega_{\Ker}=\d L_s^*\omega_{\Ker}.
\]
Moreover, left and right multiplication 
are commutative, therefore by developing the calculations we obtain
\begin{align*}
c(\xi)\cdot \omega_{\Ker}\left(s^{-1}\xi\ell(x)^{-1}\right)&=\d L_s^*\d R_{\ell(x)}^*\omega_{\Ker}(\xi)\\
&=\d R_{\ell(x)}^*\d L_s^*\omega_{\Ker}(\xi)\\
&=c(s^{-1}\xi)\cdot \omega_{\Ker}\left(s^{-1}\xi\ell(x)^{-1}\right).
\end{align*}
Thus, $c(\xi)=c(s^{-1}\xi)$ for any $s\in G^z$ which is equivalent to
the thesis of the lemma.\fine

With a slight abuse of notation we denote by $c(x)$ the value of $c(\xi)$ on the fiber $G^x$,
obtaining a smooth and everywhere non-null function $c\colon M\to \R$.
We can resume the result above with the formula
\begin{equation}\label{forc}
\d R^*_{\ell(x)}\omega_{\Ker}=c(x)\cdot \omega_{\Ker}.
\end{equation}\newline

Formulas \eqref{eqE} and \eqref{forc} allow to prove Theorem \ref{teoimp}.
Indeed, by what we showed, if $E(\xi)=x$, then
\begin{align*}
\omega_G(\xi)&=(\omega_{\bot}\wedge\omega_{\Ker})(\xi)\\
&=\left({\rho(x)}\cdot \d_\xi E^*\vol_M\right)\wedge\left({c(x)}\cdot (\d R_{\ell(x)}^{-1})^*\omega_{\Ker}\right)\\
&=\overline\rho(x)\cdot \d_\xi \psi^*\left(\vol_M\wedge \omega_{\Ker}\right),
\end{align*}
where we used the notation
\begin{align}\label{rho}
\rho(x)&=\frac{1}{\sqrt{\det\left(\mc X(x)\mc X(x)^T\right)}}\\
\label{rhobar}
\overline \rho(x)&=\rho(x)\cdot{ c(x)}
\end{align}
and the fact that
\[
\d_\xi\psi=\left(\d_\xi E|_{\Ker(\d E)^\bot}\right)\oplus \d R_{\ell(x)}^{-1}.
\]\newline

\section{Saturating fundamental solutions}\label{sfs}

\subsection{Lifting differential operators}\label{lift}
 Consider a differential operator $\mc L$ over $M$

\[
\mc L=\sum_{|\al|\leq k} r_\al(x)X^\al,
\]
where the sequence $X_1,\dots, X_n$ of complete vector fields is, as before,
a basis of the Lie algebra $\g\subset \mfk X(M)$
and the coefficients $r_\al$ are smooth functions over $M$.
The map $E\colon G\to M$ induces a canonical lifting of $\mc L$ on $G$,
that is
\[
\widetilde{\mc L}:=\sum_{|\al|\leq k} r_\al(E(\xi)) \widetilde X^\al.
\]
For any vector field $\widetilde X\in \Lie(G)$, we define the vector fields
 $\widetilde X^K(\xi)$ and $\widetilde X^\bot(\xi)$,
its projections to $\Ker(\d_\xi E)$ and $\Ker(\d_\xi E)^\bot$ respectively.
If $X\in \g$, $\widetilde X=L^{-1}X$
and $E(\xi)=x$,
then we have
\begin{equation}\label{eq1}
\d_\xi E(\widetilde X^\bot) =X(x)
\end{equation}
as $\d E$ is the first coordinate projection in $TG\cong TM\times TG^z$.

Regarding the vector component along the $G^x$ fiber, first we observe that
the vectors $\widetilde X^K$ are invariant with respect to left multiplication by elements of $G^z$.
Indeed, by Remark \ref{rmkimp} the subspace $\Ker(\d_\xi E)\subset \Lie(G)=T_\xi G$ is invariant
with respect to left multiplication $L_s$ if $s\in G^z$. At the same time, the vectors in $\Lie(G)$ and the scalar product
that we introduced on it, are invariant with respect to left multiplication
by any $G$ element. This proves that the projection $\widetilde X^K$ of $\widetilde X$ on $\Ker(\d E)$ only
depends on the fiber $G^x$ where it is taken.\newline

As a consequence of all this, if $\widetilde X^K(\xi)$ is the vertical component
of a vector $\widetilde X\in \Lie(G)=T_\xi G$ at any point $\xi$ such that $E(\xi)=x$, then
the coordinate change described in \eqref{natid}
identifies it with
\[
\d R_{\ell(x)}^{-1}\widetilde X^K(\xi)\in \Ker(\d_{\xi\ell(x)^{-1}}E)=T_{\xi\ell(x)^{-1}}G^z\cong\Lie(G^z)\subset \Lie(G).
\]
Observe that the map $\widetilde X\mapsto \d R_{\ell(x)}^{-1}\widetilde X^K$ is linear
for any $\xi\in G$, therefore we can write it as a linear (degenerate)
morphism
\[
\mc M(\xi)\colon \Lie(G)\to \Lie(G^z),
\]
varying smoothly with respect to $\xi$.
\begin{lemma}
The map $\mc M$ introduced above is invariant along any $G^x$ fiber.
\end{lemma}
\proof
For any $\widetilde X\in \Lie(G)$, $s\in G^z$ and $x\in M$,
\begin{align*}
\mc M(\xi)\widetilde X(s\xi\ell(x)^{-1})&=\d L_s\d R_{\ell(x)}^{-1}\widetilde X^K(\xi)\\
&=\d R_{\ell(x)}^{-1}\d L_s \widetilde X^K(\xi)\\
&=\d R_{\ell(x)}^{-1}\widetilde X^K(s\xi)\\
&=\mc M(s\xi)\widetilde X(s\xi\ell(x)^{-1}),
\end{align*}
where we used the commutativity of left and right multiplication.
Since this is true for any $\widetilde X\in\Lie(G) $ and $s\in G^z$, 
then $\mc M$ is constant on any fiber $G^x$.\fine

With another abuse of notation we write $\mc M(x)$ for the map $\mc M$ on the
fiber $G^x$, obtaining a smooth and everywhere maximum rank morphism
$\mc M\colon M\to \Lin(\Lie(G),\Lie(G^z))$ such that
\begin{equation}\label{eq2}
\d R_{\ell(x)}^{-1}\widetilde X^K=\mc M(x)\widetilde X.
\end{equation}
Thanks to \eqref{eq1} and \eqref{eq2} we write down the lifting of any 
vector field $X\in \g$ in the coordinate system of $M\times G^z$.
In particular, $X$ is lifted to $\widetilde X$ on $TG$ and to
\[
(\d E \widetilde X^\bot, \d R_{\ell\circ E}^{-1}\widetilde X^K)=(X,(\mc M\circ E)\widetilde X)
\]
on $TM\times TG^z$. In particular the lifted operator $\widetilde{\mc L}$
can be represented as
\begin{align*}
\widetilde{\mc L}&=\sum_{|\al|\leq k}r_\al(x)(X+\mc M(x)\widetilde X)^\al\\
&=\sum_{|\al|\leq k}r_\al(x)X^\al+
\sum_{\substack{\beta+\gamma=\al\\ |\gamma|\geq1}}r_\al(x)X^\beta(\mc M(x)\widetilde X)^\gamma.
\end{align*}
If we consider a basis $Y_1,\dots, Y_p$ of $\Lie(G^z)$ and consider the representation
in this basis of~$\mc M$, we get in particular
$\mc M(x)\widetilde X_i=\sum_{j=1}^p\mc M_{ij}(x)Y_j$ for any $i=1,\dots m$.
After rearranging the terms, we get
\begin{align*}
\widetilde{\mc L}&=\sum_{|\al|\leq k}r_\al(x)X^\al+
\sum_{\substack{|\beta+\gamma|\leq k\\ |\gamma|\geq1}}r_{\beta,\gamma}(x)X^\beta Y^\gamma\\
&=\mc L+ R.
\end{align*}
\begin{rmk}
Observe that with this representation of the lifted operator,
if we consider the coordinates $(x,s)$ on $M\times G^z$, the $X$ vector fields
only act on the $x$ coordinates while the $Y$ vector fields only
act on the $s$ coordinates. This separation of the coordinates
will be crucial in generalizing the
 Biagi-Bonfiglioli technique.
\end{rmk}

As a consequence, if we define in the usual way the dual operator $\mc L^*$
over $M\times G^z$,  it has a similar form
\begin{equation}\label{dualo}
\mc L^*=\sum_{|\al|\leq k}r^*_\al(x)X^\al+\sum_{\substack{|\beta+\gamma|\leq k\\ |\gamma|\geq1}}r_{\beta,\gamma}^*(x)X^\beta Y^\gamma.
\end{equation}
\begin{ex}\label{exsin}
Consider the operator $\mc L=\de_{x_1}^2+\left(\sin(x_1)\de_{x_2}\right)^2$.
As in the case of the Grushin operator the Lie algebra generated by $X_1=\de_{x_1}$
and $X_2=\sin(x_1)\de_{x_2}$ verifies the H\"ormander's condition everywhere,
but in this case $\g=\Lie(X_1,X_2)=\langle X_1,X_2,X_3\rangle$ is not nilpotent.
Indeed $[X_1,X_2]=X_3=\cos(x_1)\de_{x_2}$ while $[X_1,X_3]=-X_2$.

If $G$ is the simply connected Lie group whose Lie algebra is $\g$, and we put
on $\Lie(G)$ the metric induced by the orthonormal basis $\widetilde X_1,\widetilde X_2, \widetilde X_3$,
then the map $E\colon G\to \R^2$
respects
\begin{align*}
\Ker(\d E)&=\langle X_3'=\cos(x_1)\widetilde X_2-\sin(x_1)\widetilde X_3\rangle\\
\Ker(\d E)^\bot&=\langle \widetilde X_1,\  X_2'=\sin(x_1)\widetilde X_2+\cos(x_1)\widetilde X_3\rangle.
\end{align*}
If $\xi_1,\xi_2,\xi_3$ are the exponential coordinates of the first type
associated to $\widetilde X_1,\widetilde X_2,\widetilde X_3$,
then the coordinates $x_1,x_2,x_3$ associated to the vector fields $\widetilde X_1,X_2',X_3'$ in $\Lie(G)$
are
\begin{align*}
x_1&=\xi_1\\
x_2&=\frac{\xi_2}{\xi_1}(1-\cos(\xi_1))+\frac{\xi_3}{\xi_1}\sin(\xi_1)\\
x_3&=\frac{\xi_2}{\xi_1}\sin(\xi_1)-\frac{\xi_3}{\xi_1}(1-\cos(\xi_1)).
\end{align*}
We can now observe, by recovering the group operation in these coordinates, that $G\cong \R^1\rtimes \R^2$,
the universal cover of the roto-translation group $SE(2)$. Moreover, we can write down the liftings
in the new coordinates,
\begin{align*}
L^{-1}X_1&=\de_{x_1}\\
L^{-1}X_2&=\sin(x_1)\de_{x_2}+\cos(x_1)\de_{x_3}.
\end{align*}\newline
\end{ex}

\subsection{Main result}
We are going to introduce the notion of a fundamental solution for the lifted
operator $\widetilde{\mc L}$ and state that if such a solution exists, then
it is possible to saturate it and obtain a fundamental solution for $\mc L$
on the base manifold $M$.\newline

\begin{rmk}\label{rmk_hyp}
A fundamental solution of a differential
operator $\widetilde{\mc L}$ over the group $G$ is a smooth function 
\[
\widetilde \Gamma_G\colon (G\times G)\backslash \Delta(G)\to \R,
\]
where $\Delta(G)\subset G\times G$ is the diagonal,
$\widetilde \Gamma_G$ is in $L^1_{\loc}(G\times G)$ and for any $\widetilde\varphi\in C^\infty_0(G)$
\begin{equation}\label{eqff}
\int_G\widetilde \Gamma(\xi;\eta)\widetilde{\mc L}^*\widetilde\varphi(\eta)\omega_G(\eta)=-\widetilde\varphi(\xi),
\end{equation}
where $\omega_G$ is the already introduced
left invariant Haar volume form.
Moreover $\widetilde \Gamma$ respects the following properties,
\begin{enumerate}
\item $\widetilde\Gamma(\xi;\eta)>0$ for any $\xi\neq \eta$ in $G$;
\item $\widetilde \Gamma(\xi;-)\colon G\backslash\{\xi\}\to \R$ is smooth and  $\widetilde{\mc L}$-harmonic, meaning that $\widetilde{\mc L}\widetilde \Gamma(\xi;-)=0$ over $G\backslash\{\xi\}$;
\item if we consider the restricted function 
\[
\left.\widetilde\Gamma(\xi;-)\right|_{G^y}
\]
over the fiber $G^y$ for some $y\in M$ and $y\neq E(\xi)$, then $\widetilde\Gamma(\xi;-)|_{G^y}$ is in $L^1(G^y)$;
\item for any $K\subset M$ compact subset and $E^{-1}(K)\subset G$ its preimage, then
$\widetilde \Gamma(\xi;-)$ is in $L^1(E^{-1}(K))$.
\end{enumerate}
\end{rmk}


\begin{teo}\label{main}
Consider a simply connected
Riemannian manifold $M$ with a differential operator $\mc L=\sum r_\al\cdot X^\al$
on it, and such that the Lie algebra $\g\subset\mfk X(M)$ generated by $X_1,X_2,\dots,X_q$
satisfies the hypothesis of Remark \ref{rmk1}. Moreover, let $G$ be the unique connected and simply connected
Lie group such that $\Lie(G)\cong\g$.

If a fundamental solution $\widetilde{\Gamma}(\xi;\eta)$ exists for the lifted operator
$\widetilde{\mc L}$ over $G$ with the properties listed above,
then a fundamental solution $\Gamma(x;y)$ for the differential
operator $\mc L$ over $M$ exists too, and in particular the latter is obtained as an integral
\begin{equation}\label{maineq}
\Gamma(x;y):=\rho(y)\cdot \int_{G^y}\widetilde \Gamma(\xi;\eta)\omega_{\Ker}(\eta),\ \ \forall x,y\in M,\ x\neq y,
\end{equation}
where $\rho\colon M\to \R$ is the smooth and everywhere non-null function introduced
in \eqref{rho}, while
$\xi$ is any element of~$G^x$. Moreover, the function $\Gamma$ 
 verifies the following properties,
\begin{enumerate}
\item $\Gamma(x;y)>0$ for any $x\neq y$ in $M$;
\item $\Gamma(x;-)\colon M\backslash\{x\}\to \R$ is smooth and $\mc L\Gamma(x;-)=0$ over the same domain;
\item $\Gamma(x;-)$ is in $L^1_{\loc}(M)$.
\end{enumerate}
\end{teo}
\begin{rmk}
As already said, $\omega_{\Ker}$ is the restriction
of the volume form $\omega_G$ along the vertical bundle $TG^x$ for any $x\in M$.
Therefore in the formula \eqref{maineq}
the fundamental solution $\widetilde \Gamma$ is integrated along
 the $E$-fibers.
\end{rmk}

\begin{rmk}\label{rmk_main}
Preliminary to the proof, we observe that the integral in \eqref{maineq} depends
on the chosen point $\xi\in G$. Anyway, in many cases its value
is independent on this choice. For example, if the fundamental solution
is in the form $\widetilde \Gamma(\xi,\eta)=\widetilde \Gamma_G(\xi^{-1}\eta)$ with $\widetilde\Gamma_G$ a locally
integrable function on $G$, smooth outside a pole at the identity $e$.
This is the case for the solution of the sub-Laplacian on a Carnot group.

In general, if we have some additional condition that gives the unicity of $\widetilde\Gamma$,
for example some decay property at $\infty$ when $G\cong \R^n$, then
we can prove that $\Gamma$ is the same for any $\xi$.
In this case, for any $g\in G$ we have
\begin{equation}\label{eqimp}
\widetilde \Gamma(g \xi;g \eta)=\widetilde\Gamma(\xi;\eta)\ \ \forall \xi,\eta\in G,\ \xi\neq \eta.
\end{equation}
In order to prove this equation, denote by $L_g\widetilde\varphi(\xi)$ the function $\widetilde\varphi(g\xi)$
 for any compactly supported smooth function $\widetilde\varphi\colon G\to \R$ and $g\in G$,
then
\begin{align*}
\int_G\widetilde \Gamma(g\xi ;g\eta )\widetilde{\mc L}^*\widetilde\varphi(\eta)\omega_G&= 
\int_G\widetilde \Gamma(g\xi ;g\eta )\widetilde{\mc L}^*(L_{g^{-1}}\widetilde\varphi(g\eta ))\omega_G\\
&=\int_G\widetilde \Gamma(g\xi ;\eta')\widetilde{\mc L}^*(L_{g^{-1}}\widetilde\varphi(\eta'))\omega_G\\
&=- L_{g^{-1}}\widetilde\varphi(g\xi )\\
&=-\widetilde\varphi(\xi),
\end{align*}
where we used the fact that if $\eta'=g\eta $, then $\omega_G(\eta')=\omega_G(\eta)$
by left invariance.
 By the unicity of $\widetilde \Gamma$,
this proves Equation~\eqref{eqimp}.

This allows to conclude that Equation \eqref{maineq} is independent of the $\xi$ choice. Indeed,
if we consider another element $\xi'$ of $G^x$, this means that $\xi'=s \xi$ with $s\in G^z$.
We observe that $\eta'=s\eta$ is again in $G^y$ if $\eta\in G^y$, therefore we use this coordinate change
and obtain
\[
\int_{G^y}\widetilde \Gamma(\xi';\eta')\omega_{\Ker}(\eta')=\int_{G^y} \widetilde\Gamma(\xi;\eta)\omega_{\Ker}(s\eta)
=\int_{G^y}\widetilde\Gamma(\xi;\eta)\omega_{\Ker}(\eta).
\]
Here we used that $\omega_{\Ker}(\eta')=\omega_{\Ker}(\eta)$ by left invariance
with respect to $G^z$ multiplication (see proof of Lemma \ref{lemc}).\newline
\end{rmk}

\subsection{The saturation method}
In this section and the following we are going to prove Theorem \ref{main}
via a saturation method that follows the idea introduced by Biagi and Bonfiglioli in~\cite{biabon17}.

We will use the notations $(x,s)$ or $(y,t)$ for the coordinates 
on $M\times G^z$.
In particular we choose over $G^z$ the exponential coordinates ($s$ or $t$)
with respect a basis of $\Lie(G^z)$. Therefore, for the volume form
over any fiber $G^x$ we can write $\omega_{\Ker}=\d s$ or $\d t$
in order to make the calculations clearer.

With another minor abuse of notation we denote by $\widetilde \Gamma(x,s;y,t)$
the ``read'' of the fundamental solution $\widetilde\Gamma$ in the latter
coordinate system.
Observe that the $L^1$ hypothesis in Theorem~\ref{main}
translates as the fact that $\widetilde \Gamma(x,s;y,-)\colon G^z\to \R$ is of class $L^1(G^z)$ for any $(x,s)$ and $y\neq x$,
while $\widetilde \Gamma(x,s;-,-)$ is of class $L^1(K\times G^z)$ for any
compact subset $K\subset M$.
Through the coordinate change formula \eqref{forc}
we have
\begin{align*}
\int_{G^y}\widetilde\Gamma(\xi;\eta)\omega_{\Ker}(\eta)&=
\int_{G^z}\widetilde \Gamma(x,s;y,t)\omega_{\Ker}(R_{\ell(y)}(t))\\
&=c(y)\cdot \int_{G^z}\widetilde\Gamma(x,s;y,t)\d t,
\end{align*}
and this allows to rewrite Equation~\eqref{maineq},
\begin{equation}\label{maineq2}
\Gamma(x;y)={\overline\rho(y)}\cdot \int_{G^z}\widetilde \Gamma(x,s;y,t)\d t,\ \ \forall x,y\in M,\ x\neq y.
\end{equation}
where we used the definition \eqref{rhobar}.\newline

We write also the defining equation of the fundamental solution $\widetilde \Gamma$ in the new coordinate setting.
For any function $\widetilde\varphi\in C^\infty_0(G)$ we use the same notation by writing $\widetilde \varphi(x,s)=\widetilde \varphi(\xi)$
if $(x,s)=\psi(\xi)$. Therefore,
\begin{align*}
-\widetilde\varphi(x,s)=-\widetilde\varphi(\xi)&=\int_{G}\widetilde\Gamma(\xi;\eta)\widetilde{\mc L}^*\widetilde\varphi(\eta)\omega_G\\
&=\int_{G}\widetilde\Gamma(\xi;\eta)\widetilde{\mc L}^*\widetilde \varphi(\eta)\overline\rho(E(\eta))\d \psi^*(\vol_M\wedge \omega_{\Ker})\\
&=\int_{M\times G^z}\widetilde \Gamma(x,s;y,t)\widetilde{\mc L}^*\widetilde \varphi(y,t)\overline\rho(y)(\vol_M\wedge \d t),
\end{align*}
where we used the result of Theorem \ref{teoimp}.\newline

We can now prove Theorem \ref{main}.
In order to have a simpler notation, in the following we impose $s=e$ in the above formula,
which means imposing $\psi(\xi)=(x,e)$, while we will use the notation $\psi(\eta)=(y,t)$.
The proof is the same for any other choice of $s$, and in many cases also the value of the integral,
as we exposed in Remark \ref{rmk_main}.

Consider two compactly supported functions $\varphi\in C^\infty_0(M)$ and $\theta\in C^\infty_0(G^z)$
such that $\theta(e)=1$,
and define the
 product function
\[
\widetilde \varphi:=\varphi \cdot \theta.
\]
Therefore, by definition of fundamental solution and by Equation \eqref{dualo},
\begin{align*}
\int_{G}\widetilde\Gamma(\xi;\eta)\cdot \widetilde{\mc L}^*(\widetilde \varphi(\eta))\omega_G&=
 -\widetilde\varphi(\xi)= -\varphi(x)\cdot \theta(e)=-\varphi(x)=\\
  &=\int_{M\times G^z}\widetilde \Gamma\cdot\theta\cdot (\mc L^*\varphi)\cdot  \overline\rho\cdot (\vol_M\wedge \d t)
  +\int_{M\times G^z}\widetilde \Gamma\cdot R^*(\varphi\theta)\cdot \overline \rho\cdot (\vol_M\wedge \d t),
  \end{align*}
We denote by $I$ and $II$ the two integrals in the right side, 
and we want to prove that for an opportune sequence of functions $\theta_j$ with $j\to+\infty$,
we have two convergence results
\begin{align*}
I&\to \int_M\Gamma(x;y) \mc L^*\varphi(y)\vol_M\\
II&\to 0.
\end{align*}\newline

\subsubsection{The convergence of $I$}
We suppose that the $\theta_j\colon G^z\to \R$ are compactly supported cut-off functions
such that $\{\theta_j=1\}\uparrow_j G^z$ and $\theta_j(e)=1$ for each $j$.
Observe that if
$\supp(\varphi)\subset K$ a compact subset of $M$, then we have the inequality
\[
\left|\mc L^*\varphi\cdot \widetilde\Gamma\cdot \theta_j\cdot\overline\rho\right|\leq C\widetilde \Gamma,
\]
where $C$ is a constant depending on $K$.
As $\widetilde\Gamma(x,e;-,-)$ is integrable over $K\times G^z$, 
by dominated convergence we have
\[
I_j\xrightarrow{j\to+\infty}\int_{M\times G^z}\widetilde \Gamma(x,e;y,t)\cdot \mc L^*\varphi(y)\cdot\overline \rho(y)(\vol_M\wedge \d t).
\]
By the integrability of $\widetilde\Gamma(x,e;y,-)$ over $G^z$ and the Fubini's Theorem,
we get an equivalent formulation of the same result,
\[
I_j\to \int_M\mc L^*\varphi(y)\cdot\left(\int_{G^z}\widetilde \Gamma(x,e;y,t)\d t\right)\overline\rho(y)\vol_M
=\int_M\Gamma(x;y)\cdot \mc L^*\varphi(y)\vol_M
\]
as we intended to prove.\newline

\subsubsection{The convergence of $II$}\label{c2}
As showed in \eqref{dualo},
\[
R^*(\varphi(y)\cdot \theta_j(t))=\sum_{|\gamma|\geq1}r^*_{\beta,\gamma}(y)\cdot X^\beta\varphi(y)\cdot  Y^\gamma\theta_j(t).
\]
If the $\theta_j\colon G^z\to \R$  are compactly supported cut-off 
functions invading $G^z$, then
$R^*(\varphi\theta_j)$ tends pointwise to $0$ over $G^z$ because
the term $Y^\gamma\theta_j$ tends pointwise to $0$.
Therefore, if we can use again the dominated convergence argument,
we can conclude. For this reason
we want to prove that $|\overline\rho(y)\cdot r^*_{\beta,\gamma}(y)\cdot X^\beta\varphi\cdot Y^\gamma\theta_j|$
is bounded over $K\times G^z$. 
In fact, it suffices to bound the terms $| Y^\gamma\theta_j|$
because the other terms depend  only on the coordinate $y$ and
must be bounded over the compact $K$.

Therefore, if we build a sequence $\theta_j$ such that $|Y^\gamma \theta_j|\leq C$ for any $j$ and for any $\gamma$ of bounded length ($C$
depending on the compact set $K$), then
\[
\left|\widetilde \Gamma\cdot \overline\rho\cdot R^*(\varphi \theta_j)\right|\leq\widetilde\Gamma\cdot \overline\rho(y)\cdot
\sum_{|\gamma|\geq1}\left|r_{\beta,\gamma}^*(y)X^\beta\varphi(y)\right|\cdot|Y^\gamma\theta_j(t)|\leq C'\cdot \widetilde\Gamma
\]
for some constant $C'$. As $\widetilde\Gamma(x,e;-,-)$ is in $L^1(K\times G^z)$ by hypothesis,
 we can conclude
 \[
 II_j\to0.
 \]\newline

\subsection{The $\theta$ sequence}
In this section we introduce a ``good'' sequence of functions $\theta_j$
in order to complete the previous proof. We recall that we want a sequence $\theta_j\colon G^z\to \R$
for $j=1,2,\dots$ such that
\[
\theta_j(e)=1\ \forall j\ \ \m{and}\ \ \{\theta_j=1\}\uparrow_j G^z.
\]
Moreover we want the terms $ Y^\gamma\theta_j$ to be bounded,
uniformly in $j$, for $\gamma$ of bounded length.

Let's introduce a cut-off function $\theta_0 \colon G^z\to \R$ such that
$\supp \theta_0\subset K$
where $K$ is a compact $G$ subset, and there exists an open neighborhood 
$U\subset K$ of $e$ such that $\theta_0|_U\equiv 1$.
For any $g\in G^z$ we define $\theta_g:=g_*\theta_0=\theta_0(g^{-1}\cdot-)$ and we denote
by $U_g, K_g$ the pushforward of the $U$ and $K$ set respectively,
\[
U_g=g_*U=\{x|\ g^{-1}x\in U\},\ \ K_g=g_*K.
\]
By construction the $K_g$ are all compact.

Consider a vector field $Y\in \Lie(G^z)$ seen as a left invariant
vector field in $TG^z$, we observe that for any $g,h\in G^z$,
\begin{equation}\label{dg}
Y_h\theta_g=(\d L_{h} Y_e) \theta_0(g^{-1}\cdot-)=Y_e\theta_0(g^{-1}h\cdot-)=Y_{g^{-1}h}\theta_0.
\end{equation}
\begin{rmk}\label{bounds}
This implies that the bounds of the function $Y^\gamma\theta_g\colon G^z\to \R$ for
some $\gamma$ multi-index, are the same
of the function $Y^\gamma\theta_0\colon G^z\to \R$ for any $g\in G^z$, 
and both are clearly bounded because their support are included in $K_g$ and $K$.
In this work, we will only consider the case of $\gamma=0$,
but this result gives a tool for developing estimates in the $M\times G^z$
coordinate system.\newline
\end{rmk}

Consider a basis $Y_1,\dots, Y_p$ of $\Lie(G^z)$,
and the metric that makes these vectors orthonormal. For
$\varepsilon$ sufficiently small we  identify a small ball of radius $\varepsilon$ around $0\in \Lie(G^z)$, with
a neighborhood of any $g\in G^z$.
We denote by $B(g,\varepsilon)\subset G^z$ this ball-neighborhood.
We define,
\[
U^{(\varepsilon)}:=\{g\in U|\ B(g,\varepsilon)\subset U\}.
\]
We choose $\varepsilon$ sufficiently small that
$U^{(\varepsilon)}$ is an $e$ neighborhood. Observe that
as a consequence for any $h\in G^z$, 
$U_h^{(\varepsilon)}=h_*U^{(\varepsilon)}$ is a $h^{-1}$ neighborhood
and
\[
U_h^{(\varepsilon)}=\{g\in U_h|\ B(g,\varepsilon)\subset U_h\}.
\]

Consider a numerable sequence $h_1,h_2,\dots$ such that
$\left(U_{h_i}^{(\varepsilon)}\right)_i$ is a covering of $G^z$. 
Such a sequence exists because $G^z$ is second countable,
and therefore by a theorem of Lindel\"of (see \cite[Theorem VIII.6.3]{dugu66})
every open covering has a countable subcovering.\newline

We introduce
\[
\overline\theta_j:=\max_{i\leq j}\theta_{h_i}.
\]
As a consequence of Remark \ref{bounds}, the functions
$Y^\gamma\overline\theta_j$, for $\gamma$ of bounded length,
have the same bounds. Observe however that they are not defined
in every point of $G^z$ as the $\max$ might not be derivable everywhere.

In order
to obtain some smooth functions $\theta_j$ such that all the $Y^\gamma\theta_j$ have the same bounds, we consider a mollifier
$\tau_\varepsilon$ defined in a neighborhood of $0\in \Lie(G^z)$
and therefore also in a neighborhood of $e\in G^z$.
As a mollifier, $\tau_\varepsilon$ is a smooth function
verifying the following two properties
\begin{itemize}
\item $\supp\tau_\varepsilon\subset B(e,\varepsilon)$;
\item $\int_{G^z}\tau_\varepsilon=1$.
\end{itemize}
We introduce the sequence $\theta_j$ via a mollification,
that is a convolution operation with the mollifier $\tau_\varepsilon$,
\begin{equation}\label{molli}
\theta_j:=\overline\theta_j\star\tau_\varepsilon\ \ \forall j=1,2,\dots
\end{equation}
The new $\theta_j$ functions are smooth because obtained via a convolution with a smooth function.
Observe that by construction $B(e,\varepsilon)$ has finite measure
 and the functions $Y^\gamma \tau_\varepsilon$ are all bounded (for $\gamma$ of bounded length)
because the support of $\tau_\varepsilon$ is included in a compact subset.
We write
\[
|Y^\gamma \tau_\varepsilon|\leq C,\ \ |\overline \theta_j|\leq \sup |\theta_0|=C'.
\]
Therefore the $Y^\gamma\theta_j$ are bounded. Indeed,
\begin{align*}
|Y^\gamma\theta_j(h)|&=\left |\int_{G^z}Y^\gamma\tau_\varepsilon(g^{-1}h)\cdot \overline\theta_j(g)\d g\right|\\
&\leq \int_{G^z}\left| Y^\gamma\tau_\varepsilon(g^{-1}h)\right|\cdot \left|\overline \theta_j(g)\right|\d g\\
&\leq \left|B(e,\varepsilon)\right|\cdot C C',
\end{align*}
for any $\gamma$ of length $\leq k$ and for any $j$.

We also observe that 
for any $j$, we have
\[
\{\theta_j=1\}\supset \bigcup_{i\leq j}U^{(\varepsilon)}_{h_i}.
\]
By definition of the sequence $h_i$, this implies  that
the sets $\{\theta_j=1\}$ progressively invade $G^z$. Therefore we have completed 
the verification of the hypothesis on the $\theta_j$ sequence
detailed in Section \ref{c2}, and this in turn completes
the proof of the equality
\[
\int_M\Gamma(x;y)\mc L^*\varphi(y)\vol_M=-\varphi(x)\ \ \forall \varphi\in C^\infty_0(M).
\]\newline

It remains to prove that $\Gamma(x;-)$ is locally integrable
on $M$ for any $x\in M$, and the same function
is smooth and $\mc L$-harmonic on $M\backslash\{x\}$
for any $x\in M$.

The first property is a direct consequence of $\widetilde \Gamma(\xi;-)$ being in $L^1(K\times G)$ for any
compact subset $K\subset M$.

For the other properties, observe that
$\mc L\Gamma(x;-)=0$ on $M\backslash\{0\}$ as distributions, therefore
as $\mc L$ is an H\"ormander's operator it is $C^\infty$-hypoelliptic and $\Gamma(x;-)\in C^\infty(M\backslash\{x\})$.
This also implies that the equality $\mc L\Gamma(x;-)=0$ is true on $M\backslash\{x\}$ as functions.\newline

\section{The case of some non-simply connected varieties and groups}\label{nonsimply}
In this section we show an analogous result for a class of non-simply connected (but orientable) smooth manifolds.
Given a smooth manifold $N$ consider the Lie algebra
generated by smooth vector fields $\Lie(X_1,\dots,X_q)=\g\subset\mfk X(N)$,
respecting the hypothesis of Remark~\ref{rmk1}.
If $G$  is the unique
connected and simply connected Lie group such that 
$\Lie(G)\cong \g$, we denote by $L$ this isomorphism and
moreover we consider the universal cover $M$ of $N$,
therefore a regular surjection $M\to N$ with discrete fibers exists,
and $M$ is simply connected.

The vector fields $X_i$ extend naturally to smooth vector fields over $M$, so we can see $\g$ as a finite dimensional
sub-algebra of $\mfk X(M)$. If we denote by
$\mu_L$ the right $G$-action induced on $M$ by the infinitesimal generator $L\colon \Lie(G)\to \g$,
then we make the hypothesis that a discrete subgroup $H\subset G$ exists, such that
\begin{equation}
N=M\slash H.
\end{equation}
Moreover, if we denote by $R_h\colon M\to M$ the right action by $h\in H$ on $M$, then by the previous hypothesis,
\begin{equation}\label{keyeq}
\d R_hX=X,\ \ \forall h\in H,\ \forall X\in \g.
\end{equation}

\begin{teo}\label{teoact}
If $G$ is the group described above and $H$ a discrete subgroup
respecting the property \eqref{keyeq}, then
$H$ lies in the center of $G$ and therefore the right $G$-action $\mu_L$ on $M$
induces a right $\frac{G}{H}$-action on the smooth variety $N=M\slash H$.
\end{teo}
In order to prove it, we need two lemmas.

\begin{lemma}\label{lemmm1}
There is no element of the group $G$ acting trivially on $M$.
Equivalently, there exists no $g\in G$ such that $\mu_L(x,g)=x$ for every $x\in M$.
\end{lemma}
\proof
If $x_1,\dots,x_k$ are $M$ points, we denote by $\h(x_1,\dots, x_k)$ the following subspace of $\g$
\[
\h(x_1,\dots, x_k):=\left\{X\in \g|\ X(x_i)=0\ \forall i=1,\dots,k\right\}.
\]
If $n$ is the $\g$ dimension as a (real) vector space, which is the same
dimension of $G$ as a (real) manifold, and $m$ the $M$ dimension,
then there exists $n-m+1$ points $x_1,\dots,x_{n-m+1}\in M$
such that 
\begin{equation}\label{eqh}
\h(x_1,\dots x_{n-m+1})=\{0\}.
\end{equation}
Indeed, $\dim \h(x_1)=n-m$ because of the hypothesis on $\g$,
then for any $i\geq 2$ we consider $X\in \h(x_1,\dots,x_{i-1})\backslash\{0\}$ and chose
a point $x_i$ such that $X(x_i)\neq0$. It exists because $X\neq0$ as a vector field.
Therefore $\dim \h(x_1,\dots,x_{i})\leq \dim\h(x_1,\dots,x_{i-1})-1$ and by induction
we prove~\eqref{eqh}.\newline

We consider the product manifold $\widetilde M=M^{\times( n-m+1)}$ of $n-m+1$ copies
of $M$. There exists a natural right $G$-action on $\widetilde M$.
We denote by $\widetilde G\subset G$ the subgroup of elements acting trivially on the point $(x_1,\dots,x_{n-m+1})\in\widetilde M$.
As $\widetilde M$ is simply connected, the group $\widetilde G$ must be connected.
At the same time, the tangent space to $\widetilde G$ at the origin is
$\Lie(\widetilde G)=\h(x_1,\dots,x_{n-m+1})=\{0\}$.
Therefore $\widetilde G$ is trivial and the proof is concluded.\fine

\begin{lemma}\label{lemmm}
If $h$ is an $H$ element and $g$ any $G$ element, then
the right actions of $hg$ and $gh$ coincide, meaning that
\[
\mu_L(x,gh)=\mu_L(x,hg)\ \ \ \forall x\in M.
\]
\end{lemma}
\proof We fix $x\in M$ and $h\in H$ and
we prove that $\mu_L(x,g)=\mu_L(x,hgh^{-1})$.
Observe that $\nu(x,g)=\mu_L(x,hgh^{-1})$
is in fact a right $G$-action and
\[
\d\nu^{(x)}\widetilde X=\d R_h\d \mu_L^{(x)}\d L_{h^{-1}}\widetilde X,
\]
for any $\widetilde X\in \Lie(G)$.
Because of property \eqref{keyeq} and the left invariance of  any $\widetilde X$
in $\Lie(G)$,
we obtain
\[
\d \nu^{(x)}=\d \mu_L^{(x)},
\]
and because of the Fundamental Theorem \ref{fund},
this implies that the two actions are the same.\fine

As a consequence, for any $h\in H$ and $g\in G$, the commutator $ghg^{-1}h^{-1}$
acts trivially on any $x\in M$, and therefore by Lemma \ref{lemmm1},
$ghg^{-1}h^{-1}=e$. This proves that $H$ is a subgroup of the $G$ center.
Moreover, we have this further result.

\begin{lemma}\label{lemmul}
If an element $h\in H$ acts trivially on some $x\in M$, then $h=e$. More generally,
the intersection between the $G$-center and the $x$-stabilizer, is the trivial group $\{e\}$ for any $x$.
\end{lemma}
\proof Observe that
if $h$ is in the $G$-center (which is the case for any element in $H$) and $x\cdot h=x$, 
then $x\cdot g\cdot h=x\cdot h\cdot g=x \cdot g$. Therefore, as the $G$-action is transitive,
any point on $M$ is in the form $x\cdot g$ and $h$ acts trivially on $M$. By Lemma \ref{lemmm1}
this implies $h=e$.\fine
 
Consider any differential operator $\mc L_{N}$ over
the smooth variety $N=M\slash H$
in the form $\mc L_N=\sum_{|\al|\leq k} r_\al\cdot X^\al$
where the $r_\al$ are $C^\infty(N)$ coefficients. We can canonically
extend the coefficients over $M$, and therefore the whole
differential operator. We denote the extended operator by $\mc L_M$.
As in the first part of our work, this lifts to $G$, as the 
differential operator $\widetilde{\mc L}_G=\sum (r_{\al}\circ E)\widetilde X^\al$.

From now on, with a little abuse of notation 
we denote
by $E\colon G\slash H\to N= M\slash H$ the induced map on the quotients,
and by $(G\slash H)^x=E^{-1}(x)$ the preimage of any $x\in N$.
By our construction, the operatore $\widetilde{\mc L}_G$
induces naturally
a differential operator $\widetilde{\mc L}_{G\slash H}$ over $G\slash H$
which is $E$ related to $\mc L_{N}$.
\begin{ex}
We consider the analogous operator to Example \ref{exsin} defined
over $S^1\times \R$, meaning the operator $\mc L=\de_{\theta}^2+\left(\sin(\theta)\de_{x}\right)^2$.
In this case $X_1=\de_\theta$ and $X_2=\sin(\theta)\de_x$ are naturally associated to the vector fields
defined on the previous example with the same name. We have the same Lie algebra $\g=\Lie(X_1,X_2)=\langle X_1,X_2,X_3\rangle$
and the associated Lie group $G=\R\rtimes \R^2$.

The universal cover of $N=S^1\times \R$ is $M=\R^2$, and we have $M=N\slash H$
where $H$ is the discrete group of rotations $\{(2\pi k,0,0),\ k\in \Z\}\subset G$.

All our construction pass through the quotient by $H$ on $G$, and we obtain
the associated vector fields on the group $G\slash H=S^1\rtimes \R^2$, that is 
the roto-traslation group $SE(2)$,
\begin{align*}
L^{-1}X_1&=\de_\theta\\
L^{-1}X_2&=\sin(\theta)\de_x+\cos(\theta)\de_y.
\end{align*}\newline
\end{ex}

We suppose that there exists a fundamental solution
$\widetilde \Gamma_{G\slash H}$ for $\widetilde{\mc L}_{G\slash H}$
defined over the variety $((G\slash H)\times (G\slash H))\backslash \Delta$ and
such that
\begin{enumerate}
\item $\widetilde \Gamma_{G\slash H}(\xi;\eta)>0$ for any $\xi\neq \eta$ in $G\slash H$;
\item  $\widetilde \Gamma_{G\slash H}(\xi;-)\colon (G\slash H)\backslash\{\xi\}\to \R$ is $\widetilde{\mc L}_{G\slash H}$-harmonic;
\item the restricted function $\left.\widetilde \Gamma_{G\slash H}(\xi;-)\right|_{(G\slash H)^y}$ for some $y\in N$ such that $y\neq E(\xi)$
is in $L^1((G\slash H)^y)$;
\item for any $K\subset N$ compact subset and $E^{-1}(K)\subset G\slash H$ its preimage, $\widetilde \Gamma_{G\slash H}(\xi;-)$
is in $L^1(E^{-1}(K))$.
\end{enumerate}
\begin{teo}\label{main2}
Consider a simply connected
Riemannian manifold $M$.
Suppose that the Lie algebra $\g\subset\mfk X(M)$ generated by the vector fields $X_1,X_2,\dots,X_q$
satisfies the hypothesis of Remark \ref{rmk1}. Moreover, let $G$ be the unique connected and simply connected
Lie group such that $\Lie(G)=\g$. By construction a right $G$-action on $M$ exists. If $H\subset G$ is a discrete
subgroup such that $\d R_hX=X$ for any $h\in H$ and $X\in \g$, then by Theorem \ref{teoact}
there exists right $\frac{G}{H}$-action on the smooth variety $N=M\slash H$.

Consider a differential operator $\mc L_{N}=\sum_{|\al|\leq k}r_a\cdot X^\al$ on $N$.
If a fundamental solution $\widetilde{\Gamma}_{G\slash H}(\xi;\eta)$ exists for the lifted operator
$\widetilde{\mc L}_{G\slash H}$ over $G\slash H$ with the properties listed above,
then a fundamental solution $\Gamma_{N}(x;y)$ for the differential
operator $\mc L_{N}$ exists too, and in particular the latter is obtained as an integral
\begin{equation}\label{maineq3}
\Gamma_{N}(x;y):=\rho(x)\cdot \int_{(G\slash H)^y}\widetilde \Gamma_{G\slash H}(\xi;\eta)\omega_{\Ker}(\eta),\ \ \forall x,y\in N,\ x\neq y,
\end{equation}
where $\rho\colon N\to \R$ is a smooth and everywhere non-null function,
while $\xi$ is any element of the fiber $(G\slash H)^x$. Moreover, the function $\Gamma_{N}$ 
 verifies the following properties,
\begin{enumerate}
\item $\Gamma_{N}(x;y)>0$ for any $x\neq y$ in $N$;
\item $\Gamma_{N}(x;-)\colon N\backslash\{x\}\to \R$ is smooth and $\mc L_{N}\Gamma_{N}(x;-)=0$ over the same domain;
\item $\Gamma_{N}(x;-)$ is in $L^1_{\loc}(N)$.
\end{enumerate}
\end{teo}
\proof
Given the fundamental solution $\widetilde \Gamma_{G\slash H}$ for $\widetilde{\mc L}_{G\slash H}$, 
we can naturally define $\widetilde \Gamma_G(\xi;\eta)$
over $(G\times G)\backslash \overline \Delta$ where $\overline \Delta $ is the preimage by $G\times G\to (G\slash H)\times (G\slash H)$
of the diagonal $\Delta(G\slash H)$.
Observe that  $\widetilde \Gamma_{G}$ is a fundamental solution
for $\widetilde{\mc L}_G$ that satisfies the hypothesis
of Theorem \ref{main} with the slight adaptation
that  $\widetilde \Gamma_G(\xi;-)\in L^1(G^y)$ for any $y\in M$ outside
the $H$-orbit of $x=E(\xi)$. We denote by $\Gamma_{M}$ the fundamental
solution to $\mc L_M$ induced by $\widetilde \Gamma_G$ and Formula~\eqref{maineq}.

Our construction implies that 
\[
\Gamma_M(x;y)=\Gamma_M(xh_1;yh_2)\ \ \forall h_1,h_2\in H,\ x,y\in M\ xH\neq yH.
\]
Therefore it is possible to define $\Gamma_{N}\colon (N\times N)\backslash \Delta\to \R$.
Observe that $\Gamma_{N}$ is a fundamental solution for $\mc L_{N}$
respecting the properties in the theorem.

Observe that by Lemma \ref{lemmul}, the action of $H$ is trivial on any $E$-fiber.
Equivalently, $G^y$ is in bijection with $(G\slash H)^y$ via the map induced by the quotient,
for any $y\in N=M\slash H$. Observe moreover that the construction of the map $\mc X(x)$
(see Remark \ref{rmkimp}) is $H$-invariant, and therefore
the function $\rho\colon M\to \R$ is $H$-invariant too,
therefore by an abuse of notation we can define a smooth and everywhere
non-null function with the same name $\rho \colon N\to \R$.
We can conclude that Formula \eqref{maineq3}
gives the result value of $\Gamma_{N}$ because the quotient
pass through the integral.\fine

\bibliography{mybiblio}{}
\bibliographystyle{plain}

\end{document}

%% file: newcommands_eng.tex
\newtheorem{teo}{Theorem}[section]

\newtheorem{lemma}[teo]{Lemma}

\newtheorem*{teon}{Theorem}

\theoremstyle{definition}

\newtheorem{rmk}[teo]{Remark}
\newtheorem{ex}[teo]{Example}

\newcommand{\mb}{\mathbb}
\newcommand{\mc}{\mathcal}

\newcommand{\mfk}{\mathfrak}

\def\R{\mathbb{R}}

\def\N{\mathbb N}
\def\Z{\mathbb Z}

\def\H{\mb{H}}

\def\g{\mfk{g}}
\def\d{\mathrm{d}}

\newcommand{\al}{\alpha}

\newcommand{\m}{\mbox}
\newcommand{\gra}{\textbf}
\newcommand{\cor}{\textit}

\newcommand{\fine}{\qed\newline}

\DeclareMathOperator{\im}{Im}

\DeclareMathOperator{\Ker}{Ker}

\DeclareMathOperator{\vol}{Vol}

\DeclareMathOperator{\ev}{ev}

\DeclareMathOperator{\Lie}{Lie}
\DeclareMathOperator{\supp}{supp}

\def\h{\mfk{h}}
\newcommand{\de}{\partial}
\DeclareMathOperator{\Lin}{Lin}
\DeclareMathOperator{\loc}{loc}
\DeclareMathOperator{\Exp}{Exp}

